\documentclass[12pt,a4paper]{article}

\usepackage[
top=2.5cm,
left=3cm,
bottom=2.5cm,
right=2.5cm
]{geometry}

\usepackage{amstext}
\usepackage{graphicx}
\usepackage{caption}
\usepackage{subcaption}
\usepackage{amsmath, amsthm}
\usepackage{color}
\usepackage[T1]{fontenc}
\usepackage[utf8]{inputenc}
\usepackage{authblk}

\usepackage{multirow}

\usepackage[english]{babel}
\usepackage{hyperref}
\hypersetup{
 bookmarks=true,
 pdftitle={Delayed Duffing equation},
 pdfauthor={Si Mohamed Sah, Bernold Fiedler},
 colorlinks=true,
 linkcolor=blue,
 citecolor=blue,
 filecolor=blue,
 urlcolor=blue}
 
\usepackage{float}

\usepackage[textwidth=2cm,textsize=small,backgroundcolor=none]{todonotes}

\usepackage[pagewise]{lineno}
\definecolor{refkey}{rgb}{1,0,0}
\definecolor{labelkey}{rgb}{1,0,0}

\theoremstyle{plain}
\newtheorem{thm}{Theorem}[section]

\usepackage{chngcntr}
\counterwithin{figure}{section}
\counterwithin{equation}{section}
\numberwithin{equation}{section}
\numberwithin{figure}{section}
\title{Unbounded sequences of stable limit cycles\\
in the delayed Duffing equation: an exact analysis}
\author[1]{Si Mohamed Sah%\thanks{A.A@university.edu}
}
\author[2]{ Bernold Fiedler}
\author[3]{ B. Shayak}
\author[3,4]{ Richard H. Rand}
\affil[1]{Department of Mechanical Engineering, Technical~University~of~Denmark,~Denmark}
\affil[2]{Institut für Mathematik, Freie Universität Berlin, Germany}
\affil[3]{Theoretical and Applied Mechanics, Sibley School of Mechanical and Aerospace Engineering, Cornell University, Ithaca, New York 14853 USA}
\affil[4]{Department of Mathematics, Cornell University, Ithaca, New~York~14853~USA}

\date{version of \today}

\begin{document}
\maketitle

\thispagestyle{empty}

\begin{abstract}
\parindent=0cm
\parskip 4pt
\thispagestyle{empty}

\noindent

The delayed Duffing equation 
$\ddot{x}(t)+x(t-T)+x^3(t)=0$
is shown to possess an infinite and unbounded sequence of rapidly oscillating, asymptotically stable periodic solutions, for fixed delays such that $T^2<\tfrac{3}{2}\pi^2$.
In contrast to several previous works which involved approximate solutions, the treatment here is exact.
\end{abstract}

%%%%%%%%%%%%%%%%%%%%%%%%%%%%%%%%%%%%%%%%%%%%%%%%%%%%%%
\section{Introduction}\label{intr}

This work concerns a differential-delay equation (DDE) known as the delayed Duffing equation
\begin{equation}
 \ddot{x}(t)+ x(t-T) + x(t)^3 = 0\,,
\label{dde}
\end{equation}
where $T>0$ is the time delay. 
The existence of an infinite number of \emph{stable limit cycles}, i.e.~of asymptotically stable periodic solutions, in this DDE was first suggested in a paper by Wahi and Chatterjee \cite{Wahi}. 
Formally and to leading order, they performed the method of averaging and obtained a slow flow that predicted infinitely many stable limit cycles. 
In their DDE, the time delay was fixed at $T=1$. 
Mitra\&al \cite{Mitra} studied the same DDE with an added linear stiffness. 
By assuming an approximate solution in harmonic form $x(t)=A\sin(\omega t)$, they claimed that the system exhibits an infinite number of stable limit cycles for any value of the time delay $T$.  
In a paper by Davidow\&al \cite{Davidow}, the same claim was supported by a) harmonic balance, b) Melnikov's integral with Jacobi elliptic functions, and c) the introduction of damping.

Strictly speaking, however all  these works on the delayed Duffing equation (\ref{dde}) were restricted to small amplitudes of the limit cycles. In our work, we present an exact treatment of (\ref{dde}), in the limit of unboundedly large amplitudes. In particular, the previously studied infinite sequences of ``stable limit cycles'' lose stability, eventually, for delays $T$ such that $T^2>\tfrac{3}{2}\pi^2$.

Section \ref{nume} gives a brief account of the numerical integration method used for our simulations.
In section \ref{lift} we study exact periodic solutions $x_n(t)$ of a slightly generalized Duffing ordinary differential equation (ODE), with vanishing time delay $T=0$; see \eqref{ode}.
We show how the non-delay ODE solutions $x_n(t)$ of minimal (or fundamental) periods $p_n$ lift to exact solutions of the original delayed Duffing DDE \eqref{dde} with positive delay $T>0$, provided their minimal periods
\begin{equation}
\label{p_n}
p_n=2T/n
\end{equation}
are integer fractions of the double delay $2T$.
In particular we show how the more and more rapidly oscillating periodic solutions $x_n(t)$ develop unbounded amplitudes $A_n \nearrow \infty$, for $n \rightarrow \infty$.
In section \ref{ampl} we indicate how to determine the amplitudes $A_n$ of the lifted solutions $x_n$\,, numerically and by series expansions for $n \rightarrow \infty$.
Section \ref{stab} recalls our stability results from \cite{Fieetal}.
These mathematical results basically assert local asymptotic stability of the solutions $x_n(t)$, for any fixed positive delay $T$ such that $T^2<\tfrac{3}{2}\pi^2$ and for sufficiently large odd $n=1,3,5,\dots$. 
They also show instability, for sufficiently large even $n=2,4,6,\dots$.
For full mathematical details, including added linear stiffness, we refer to \cite{Fieetal}.
We conclude with numerical illustrations of our results, in section \ref{disc}, and a short summary \ref{conc}.

\textbf{Acknowledgment.} 
Just as the more mathematically inclined account in \cite{Fieetal}, the present work has originated at the \emph{International Conference on Structural Nonlinear Dynamics and Diagnosis 2018, in memoriam Ali Nayfeh}, at Tangier, Morocco. We are deeply indebted to Mohamed Belhaq, Abderrahim Azouani, to all organizers, and to all helpers of this outstanding conference series. They indeed keep providing a unique platform of inspiration and highest level scientific exchange, over so many years, to the benefit of all participants. This work was partially supported by DFG/Germany through SFB 910 project A4. Authors RHR, BS and SMS gratefully acknowledge support by the National Science Foundation under grant number CMMI-1634664.

%%%%%%%%%%%%%%%%%%%%%%%%%%%%%%%%%%%%%%%%%%%%%%%%%%%%%%

\section{Numerical integration}\label{nume}

For zero delays, $T=0$, the delayed Duffing DDE (\ref{dde}) reduces to a non-delayed ordinary differential equation (ODE) known as the classical Duffing equation.
The equation is conservative and hence exhibits a continuum of periodic orbits, rather than any asymptotically stable limit cycles. 

Even for arbitrarily small fixed positive delays, $T>0$, in contrast, approximate analysis and numerical simulations suggest that an infinite number of stable limit cycles may coexist, their amplitudes going to infinity \cite{Davidow}. 

\begin{figure}
	\centering
	\includegraphics[width=0.9\textwidth]{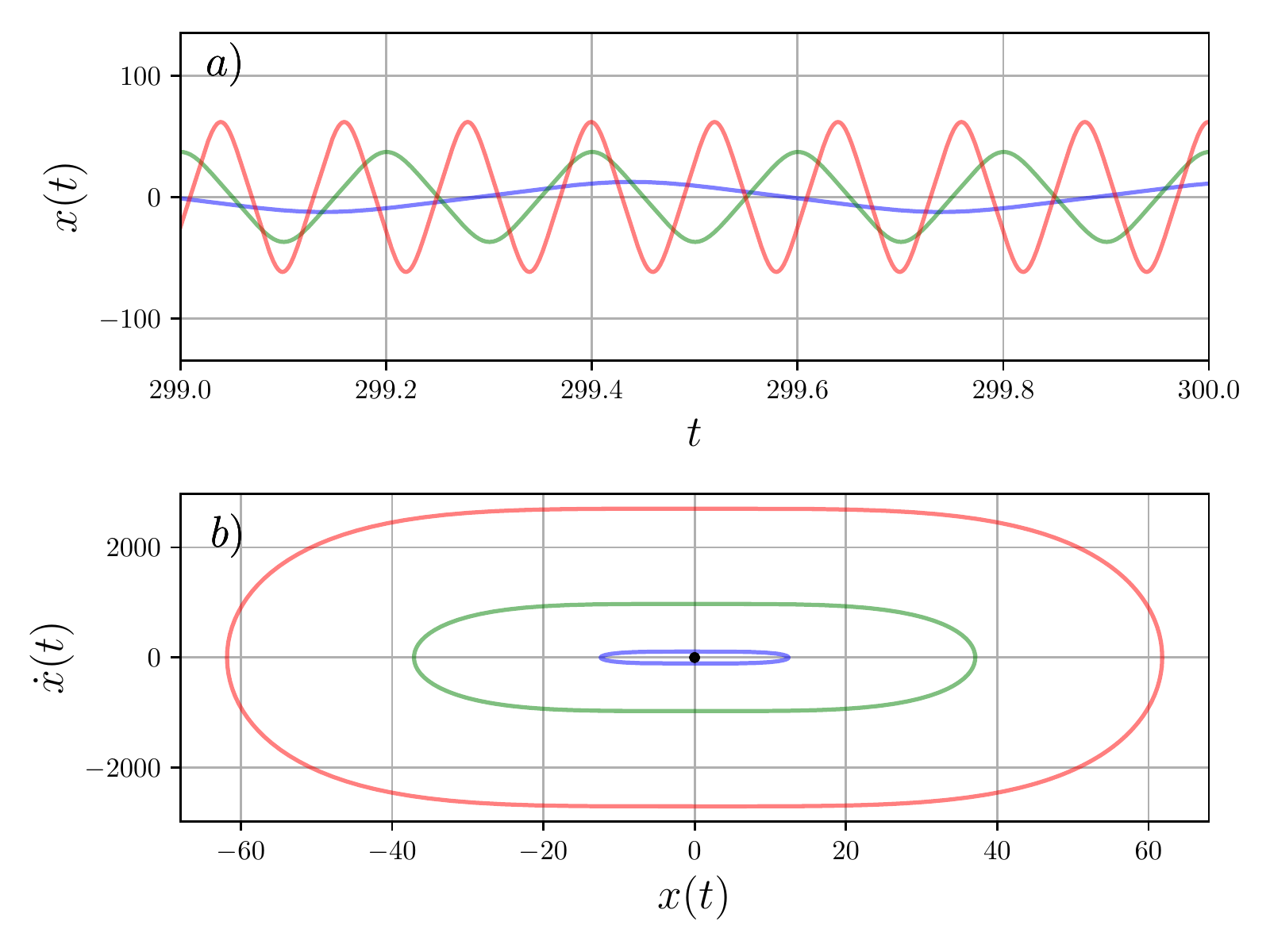} 
	\caption{(a) Time histories of some periodic solutions $x_n(t)$ for the delayed Duffing DDE (\ref{dde}) with fixed delay $T=0.3$. (b) Nested phase plane plots $(x_n(t), \dot{x}_n(t))$ of the periodic orbits $x_n$ with minimal period $2T/n$, $n=1, 3, 5$. Black dot corresponds to equilibrium point.}
	\label{fig:01}
\end{figure}

Figure \ref{fig:01} shows the time history (a) and phase plane (b) of the first three stable limit cycles obtained by numerical integration of the delayed Duffing DDE (\ref{dde}), for $T=0.3$ and with different initial conditions. 
The numerical integrations in the present work were performed using the Python library \texttt{pydelay} for DDEs \cite{{Flunkert}}. The integrator is based on the Bogacki-Shampine method \cite{Bogacki}. 
The maximal step size used to produce the plots in the present work was fixed at $\Delta t = 10^{-4}$.  
See section \ref{disc} for further numerical examples.

%%%%%%%%%%%%%%%%%%%%%%%%%%%
%%%%%%%%%%%%%%%%%%%%%%%%%%%%%%

%%%%%%%%%%%%%%%%%%%%%%%%%%%%%%%%%%%%%%%%%%%%%%%%%%%%%%

\section{Lifting periodic solutions from the non-delayed to the delayed Duffing equation}\label{lift}

In this section we show the existence of infinitely many rapidly oscillating periodic solutions of specific periods $p$ in the delayed Duffing DDE \,(\ref{dde}). 
Our approach is based on a lift of certain periodic solutions of the ordinary non-delayed Duffing ODE \eqref{ode} below, with minimal (or, fundamental) period $p$, to periodic solutions of the delayed Duffing DDE (\ref{dde}) with time delay $T$.
We will show this remarkable fact for minimal periods $p$ which are integer fractions of the doubled delay $2T=np$\,; see claim \eqref{p_n}.
We first recall some elementary facts on the non-delayed Duffing ODE, in subsection \ref{DuffingODE}. We separately address the cases of even and odd fractions $n$ in subsections \ref{neven} and \ref{nodd}, respectively.

%%%%%%%%%%%%%%%%%%%%%%%%%%%
\begin{figure}[t!]
	\centering
	%\vspace*{-2cm}
	%\hspace*{-3.2cm} 
	\includegraphics[width=1.1\textwidth]{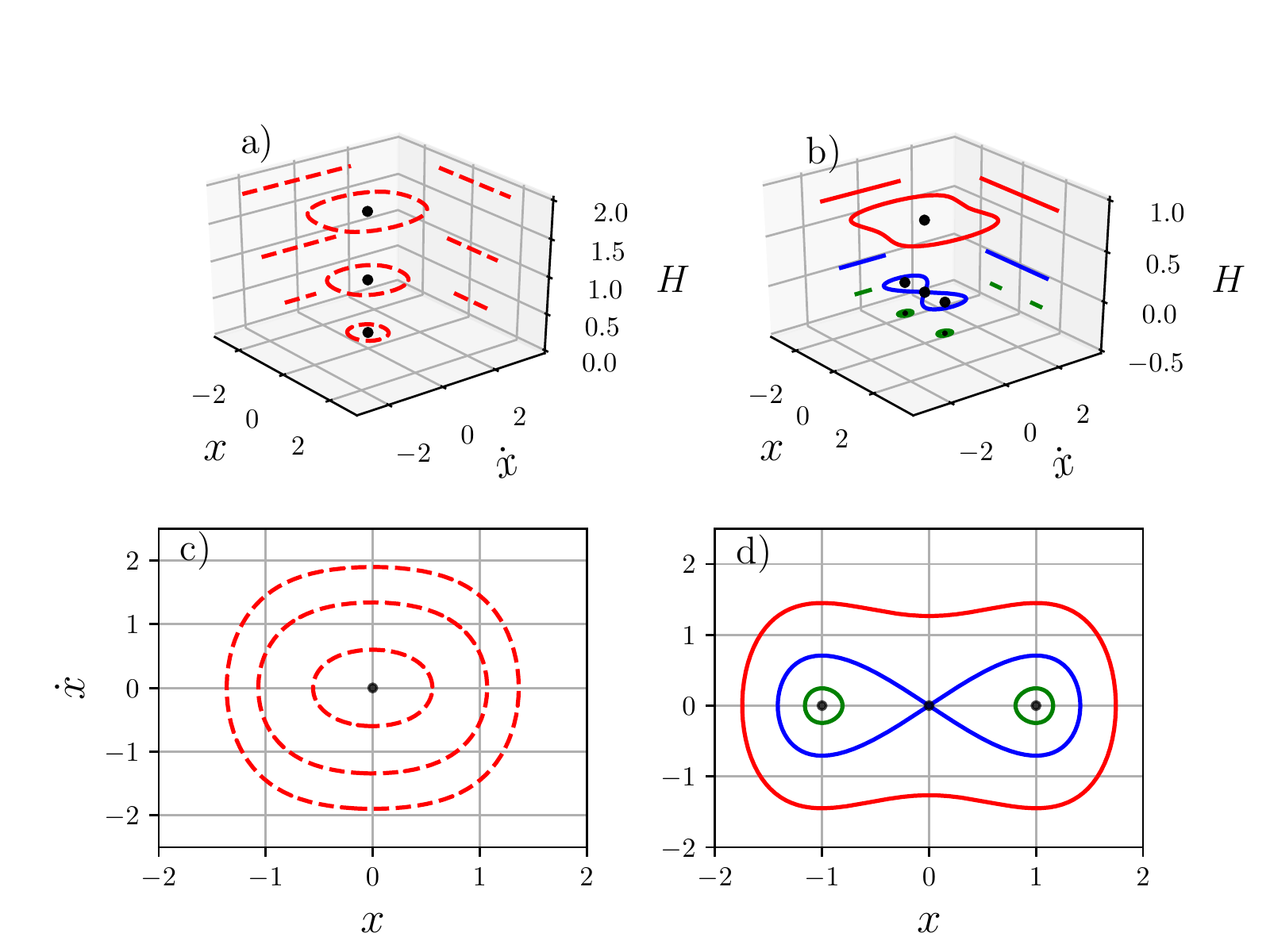}
	\caption{Three dimensional plots of Hamiltonian level sets (\ref{ham_1}) in (a,b), and projections into  the ($x, \dot{x}$) plane in (c,d), for the general non-delayed Duffing ODE (\ref{ode}). 
	(a,c) $n$ even: the single-well Duffing ODE (\ref{ode+}).
	(b,d) $n$ odd: the double-well Duffing ODE (\ref{ode-}). 
	The Hamiltonian $H$ of the double-well Duffing equation (\ref{ode-}) in (b,d) can be strictly negative (green), zero (blue), or strictly positive (red) as assumed in \eqref{H>0}. Black dots correspond to equilibrium points.  }
	\label{fig:02}
\end{figure}
%%%%%%%%%%%%%%%%%%%%%%%%%%%%%%

\subsection{General Duffing equation}\label{DuffingODE}

We consider the following two general forms of the classical Duffing ODE \cite{Kovacic}:
\begin{equation}
\ddot{x}(t) + (-1)^n x(t) + x(t)^3 = 0,  \ \,\,\,\,\,    n = 1, 2, 3, \dots ,
\label{ode}
\end{equation} 
The time-independent Hamiltonian energy of (\ref{ode}) takes the form
\begin{equation}
H(t) =\tfrac{1}{2}\,{\dot{x}^2}+ \tfrac{1}{2}\,(-1)^n\,x^2 + \tfrac{1}{4}x^4\,.
\label{ham_1}
\end{equation}
See Figure \ref{fig:02}. 
For even $n$ the Hamiltonian is always positive; see Figure \ref{fig:02}a,c. 
For odd $n$, however, see Figure \ref{fig:02}b,d: the Hamiltonian is either strictly negative (green), identically zero (blue) or strictly positive (red), depending on the ODE initial conditions.
Note how single trajectories in the $(x,\dot{x})$-plane are point symmetric to the origin, if and only if the positive energy condition
\begin{equation}
\label{H>0}
H>0
\end{equation}
is satisfied.
We assume this restriction to hold throughout our further analysis.

For $H>0$, we may time-shift solutions $(x_n(t),\dot{x}_n(t))$ of \eqref{ode} such that the initial conditions 
\begin{equation}
\label{odeic}
0<x_n(0) =: A_n\,, \qquad \dot{x}_n(0)=0, 
\end{equation}
are satisfied. 
In particular, $A_n = \max |x_n(t)|$ is the amplitude of the solution $x_n$.
For odd $n$, note how our positivity condition \eqref{H>0} requires an amplitude $A_{n}>\sqrt{2}$ in \eqref{ham_1};  see the red curve in Figure \ref{fig:02}b,d, outside the blue figure-8 shaped separatrix loops. 
The periodic closed curves fill the part of the phase space $(x,\dot{x})$ where $H>0$.
Each periodic orbit  corresponds to specific initial conditions and possesses a specific minimal period.

The exact periodic solutions of the Duffing ODE $(x_n(t),\dot{x}_n(t))$ of \eqref{ode} are easily determined.
Indeed the energy $H\equiv E$ is identically constant.
Solving \eqref{ham_1} for $\dot{x}$ and classical separation of variables therefore lead to the elliptic integrals
\begin{equation}
\label{ell}
t = \int^{x_n(t)}_{x_n(0)}\frac{dx}{\dot{x}(t)} = \pm \int^{A_n}_{x_n(t)} \frac{dx}{\sqrt{2\,E - (-1)^n\,x^2 - x^4/2 }}\,.
\end{equation}
Here we have substituted the initial condition \eqref{odeic} for $x_n(0)$.
The minimal (fundamental) period $p_{n}$ can be determined as the special case $t=p_n/4$, where symmetry implies $x_n(t)=0$:
\begin{equation}
\tfrac{1}{4}\,p_n = \int^{A_n}_{0} \frac{dx}{\sqrt{\left( 	2\,E - (-1)^n\,x^2 - x^4/2  \right)}}\,.
\label{per_n_q}
\end{equation}
Evaluating the invariant Hamiltonian $H_{n}\equiv E$ at the initial condition \eqref{odeic} provides the energy \begin{equation}
H_n = E = \tfrac{1}{2}\,(-1)^n\,A_n^2 + \tfrac{1}{4}A_n^4
\label{ham_2}
\end{equation}
and the explicit elliptic integral
\begin{equation}
\tfrac{1}{4}\,p_n = \int^{A_n}_{0} \frac{dx}{\sqrt{\left(A_n^2-x^2\right)\,\left((-1)^n + A_n^2/2 + x^2/2\right)}}\,.
\label{per_n_q_2} 
\end{equation}

\begin{figure}
	\centering
	\vspace*{-1cm}
	\includegraphics[width=.9\textwidth]{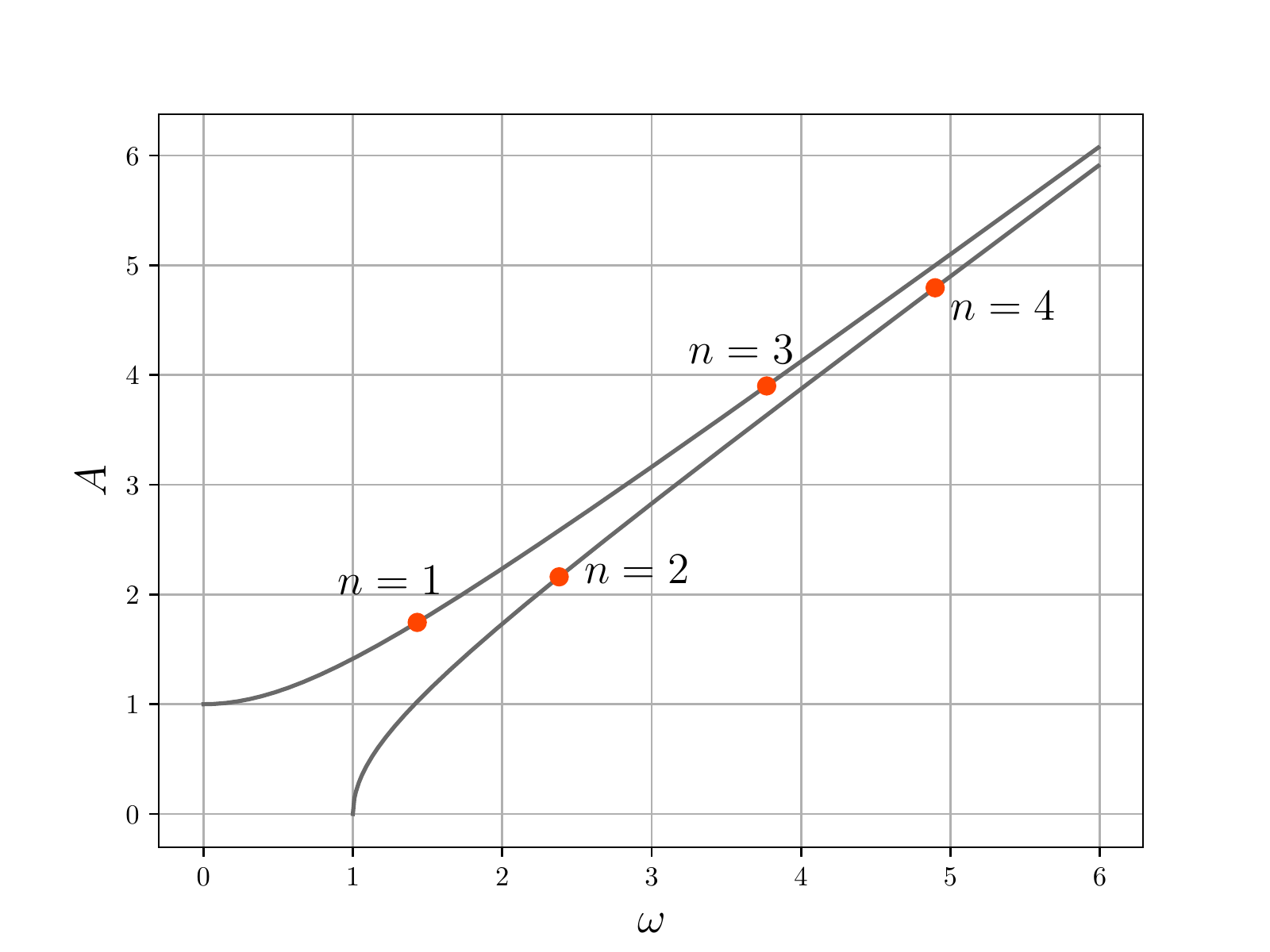} 
	\caption{  The relation between amplitude $A$ and frequency $\omega$ of the periodic solutions in the non-delayed Duffing ODE (\ref{ode}) obtained from the second equation of (\ref{init_4}). Upper curve for $n$ odd and lower curve for $n$ even. Only the marked points on these two curves correspond to periodic solutions $x_n(t)$ of the delayed Duffing DDE (\ref{dde}). The time delay for this plot is $T=3$.}
	\label{fig:03}
\end{figure}

 The elliptic integral \eqref{ell} allows us to express the exact periodic solution of the general Duffing ODE (\ref{ode}) in terms of Jacobi elliptic function as
 \begin{equation}
 x_{n}(t) = A_{n}\, \mathrm{cn}(\omega_{n}\,t, m_{n}).
 \label{sol}
 \end{equation}
  Here $\mathrm{cn}$ denotes the Jacobi elliptic cosine function. 
 The arguments $A_{n}$, $\omega_{n}$ and $0 <m_{n}< 1$ are the amplitude, the angular frequency, and the elliptic modulus, respectively. 
 The frequency $\omega_{n}$ and the modulus $m_{n}$ in the solution (\ref{sol}) are related to the amplitude $A_{n}$ such that
 \begin{equation}
 m _n= \frac{A_n^2}{2(A_n^2+(-1)^n)}~~~~~~~~~\textrm{and} ~~~~~~~~~~ \omega_n = \sqrt{A_n^2+(-1)^n}\,.
 \label{init_4}
 \end{equation} 
The minimal period \eqref{per_n_q_2} can be expressed in terms of the complete elliptic integral of the first kind $K\equiv K(m_{n})$ as
\begin{equation}
p_{n} = 4\,K/\omega_{n}.
\label{p_1}
\end{equation}
See \cite{Rand}.
Figure \ref{fig:03}  indicates the relation between amplitude and frequency for the general Duffing ODE (\ref{ode}). The two black curves are obtained from the second equation of (\ref{init_4}), and they correspond to the relation between amplitude and frequency of the periodic solutions (\ref{sol}) in the non-delayed Duffing ODE (\ref{ode}), for $n$ odd (upper curve) and $n$ even (lower curve). 
Each point represents a periodic orbit of the general Duffing ODE (\ref{ode}).
In the phase plane, each of the black curves therefore indicates a foliation by periodic solutions.
For the delayed Duffing equation \eqref{dde}, the same periodic solutions $x_{n}$ of minimal period $p_{n}=2T/n$ on the upper curve ($n$ odd)  will turn out locally asymptotically stable, for $T^2<\tfrac{3}{2}\pi^2$ and large $n$,
while large $n$ of even parity (lower curve) always turn out linearly unstable; see Theorems \ref{thmodd}, \ref{thmeven} below.

%%%%%%%%%%%%%%%%%%%%%%%%%%%
%%%%%%%%%%%%%%%%%%%%%%%%%%%%%%

Our lift construction from solutions of the non-delayed Duffing ODE (\ref{ode}) to the delayed Duffing ODE (\ref{dde}) is based on two interpretations of the mathematical expression $x(t-T)$.  
On the one hand, $x(t-T)$ represents a delay, as in (\ref{dde}).  
The same expression, on the other hand, represents a periodic solution when equated to $\pm x(t)$ by
\begin{equation}
x_{n}(t-T) = (-1)^n x_{n}(t)\,.
\label{perd}
\end{equation}
 Here $2T$ represents any (not necessarily minimal) period of the periodic solution $x(t)$. 
Indeed, any positive energy solution of the  Duffing ODE (\ref{ode}) is periodic and will automatically satisfy the periodicity condition (\ref{perd}), for some $T > 0$.  
Upon substitution of the periodicity condition (\ref{perd}), however, the non-delayed Duffing ODE (\ref{ode}) produces the delayed Duffing DDE (\ref{dde}), where now the (half) period $T$ represents the delay.  
Thus any periodic solution of the Duffing ODE (\ref{ode}) with periodicity condition \eqref{perd} lifts to a periodic solution of the DDE (\ref{dde}), for that choice of the delay $T$.  
The marked points (red) on the two black curves in Figure \ref{fig:03}, for example, correspond to periodic solutions of the delayed Duffing DDE (\ref{dde}), with delay $T=3$.
 
Actually, the non-delayed ODE Duffing equation (\ref{ode}) possesses an uncountable continuum of periodic orbits, foliating the phase plane. 
The number of periodic orbits $x_n$ which satisfy the periodicity condition \eqref{perd}, however, is (at most) countable.
In particular, our lift construction from the non-delayed Duffing ODE (\ref{ode}) to the delayed Duffing DDE (\ref{dde}) restricts the allowable points on the curves in Figure \ref{fig:03}  to a countable set and therefore produces only a countable set of periodic solutions for the delayed Duffing DDE.
We do not claim that our lift construction covers all possible periodic solutions of the DDE \eqref{dde}; in section \ref{stab} we will see indications of additional periodic solutions which cannot be obtained by our lift.

In the following we will further detail the lift construction \eqref{perd} which is based on the known exact periodic solutions (\ref{sol}) of the general Duffing ODE (\ref{ode}). We consider the two cases, $n$ even and $n$ odd, separately.

%%%%%%%%%%%%%%%%%%%%%%%%%%%
\begin{figure}[t!]
	\centering
	\vspace*{-1.2cm}
	\includegraphics[width=.95\textwidth]{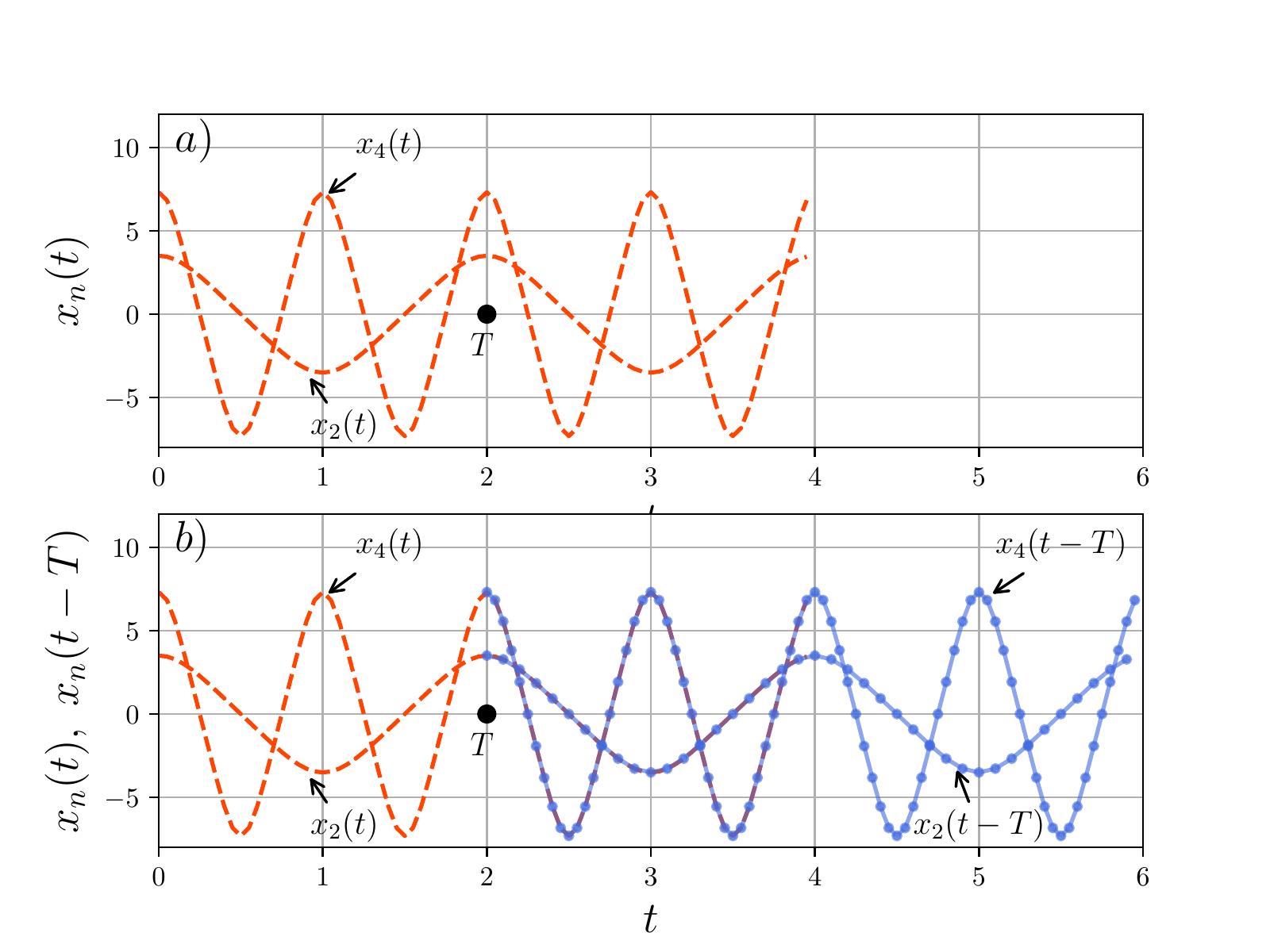}
	\caption{Solutions of the single well Duffing ODE \,(\ref{ode+}), alias even $n$ in the delayed Duffing DDE (\ref{ode}). Dashed red: solutions $x_{n}(t)$ of (\ref{ode+}). Dotted blue: shifted delayed solutions $x_{n}(t-T)$, $T=2$. }
	\label{fig:04}
\end{figure}
%%%%%%%%%%%%%%%%%%%%%%%%%%%%%%

\subsection{Even $n$}\label{neven}
For even $n$, the general Duffing ODE (\ref{ode}) reduces to the single-well case 
\begin{equation}
\ddot{x}(t) + x(t) + x(t)^3 = 0.
\label{ode+}
\end{equation}
By (\ref{sol}), the exact periodic solutions are expressed as
 \begin{equation*}
x_{n}(t) = A_{n}\, \mathrm{cn}(\omega_{n}\,t, m_{n}),
\end{equation*} 
where now (\ref{init_4}) becomes
\begin{equation}
m_{n} = \frac{A_{n}^2}{2(A_{n}^2+1)}~~~~~~~~~\textrm{and} ~~~~~~~~~~ \omega_n = \sqrt{A_{n}^2+1}.
\label{init_1}
\end{equation}
According to \eqref{p_1}, minimal periods $p$ decrease monotonically from $p=2 \pi$, at amplitude $A = 0$, to $p = 0$, for unbounded amplitudes $A \nearrow \infty$; see Figure \ref{fig:03}. 

To perform the lift from the Duffing ODE \eqref{ode+} to the Duffing DDE \eqref{dde}, we fix a time delay T (black dot $T=2$ in Figure \ref{fig:04}), a priori, such that $T < 2 \pi$. 
Then we can always find a solution (\ref{sol}) to (\ref{ode+}) with minimal period $p_{2} = T$; see solution $x_{2}(t)$ in Figure \ref{fig:04}. 
If we shift the curve of $x_{2}(t)$ to the right by $T$ we obtain a new curve  $x_{2}(t-T)$ that coincides with  $x_{2}(t) = x_{2}(t-T)$; see Figure \ref{fig:04}b.  
We can also find another solution $x_{4}(t)$, of larger amplitude $A_4>A_2$, whose minimal period is $p_{4} =T/2$. 
Shifting by $T$ we obtain a new curve $x_{4}(t-T)$ that coincides with $x_{4}(t) = x_{4}(t-T)$, see Figure \ref{fig:04}b again. 
In the  same manner, we can find infinitely many periodic solutions $x_n(t)$ with minimal periods $p_{n} = 2T/n$, for $n=2, 4, 6, \dots$. After time shift by their shared (non-minimal) period $T$ we obtain
\begin{equation}
x_{n}(t-T) = x_{n}(t),~~~~~~~~~~ \textrm{for all even }n.
\label{i_T}
\end{equation}

Substituting (\ref{i_T}) into (\ref{ode+}) lifts all those ODE Duffing solutions $x_n(t)$ to the delayed Duffing DDE (\ref{dde}), for fixed delay $T<2\pi $. 
Note how $p_n\searrow 0$ implies unbounded amplitudes $A_n \nearrow \infty$, for $n\rightarrow \infty$.
As the amplitudes $A_n$ of the periodic solutions of the Duffing ODE (\ref{ode+}) increase to infinity, the minimal periods $p_n$ decrease to zero.  
Thus we obtain an unbounded sequence of more and more rapidly oscillating periodic solutions, with minimal periods $T, T/2, T/3, \dots$, which are also periodic with (non-minimal) period $T$.
This proves our claim \eqref{p_n}, for even $n$.

\subsection{Odd $n$}\label{nodd}
For odd $n$, the general Duffing ODE (\ref{ode}) reduces to the double-well case  
\begin{equation}
\ddot{x}(t) -x(t) + x(t)^3 = 0.
\label{ode-}
\end{equation}
Any solution conserves the Hamiltonian energy
\begin{equation}
H =\tfrac{1}{2}\,{\dot{x}^2}- \tfrac{1}{2}\,x^2 + \tfrac{1}{4}x^4\,.
\end{equation}

We recall how the phase portrait of the double-well Duffing ODE (\ref{ode-}) is characterized by a figure-8 shaped separatrix $H=0$; see the blue curve in Figure \ref{fig:02}b,d. 
For positive energy $H>0$, the (red) solutions of (\ref{ode-}) oscillate around the exterior of the separatrix. 
Again, minimal periods $p$ decrease monotonically: this time from $p=\infty$, at the separatrix amplitude $A = \sqrt{2}$, to $p = 0$, for $A \nearrow \infty$.

%%%%%%%%%%%%%%%%%%%%%%%%%%%
\begin{figure}[t!]
	\centering
	\vspace*{-0.8cm}
	\includegraphics[width=.95\textwidth]{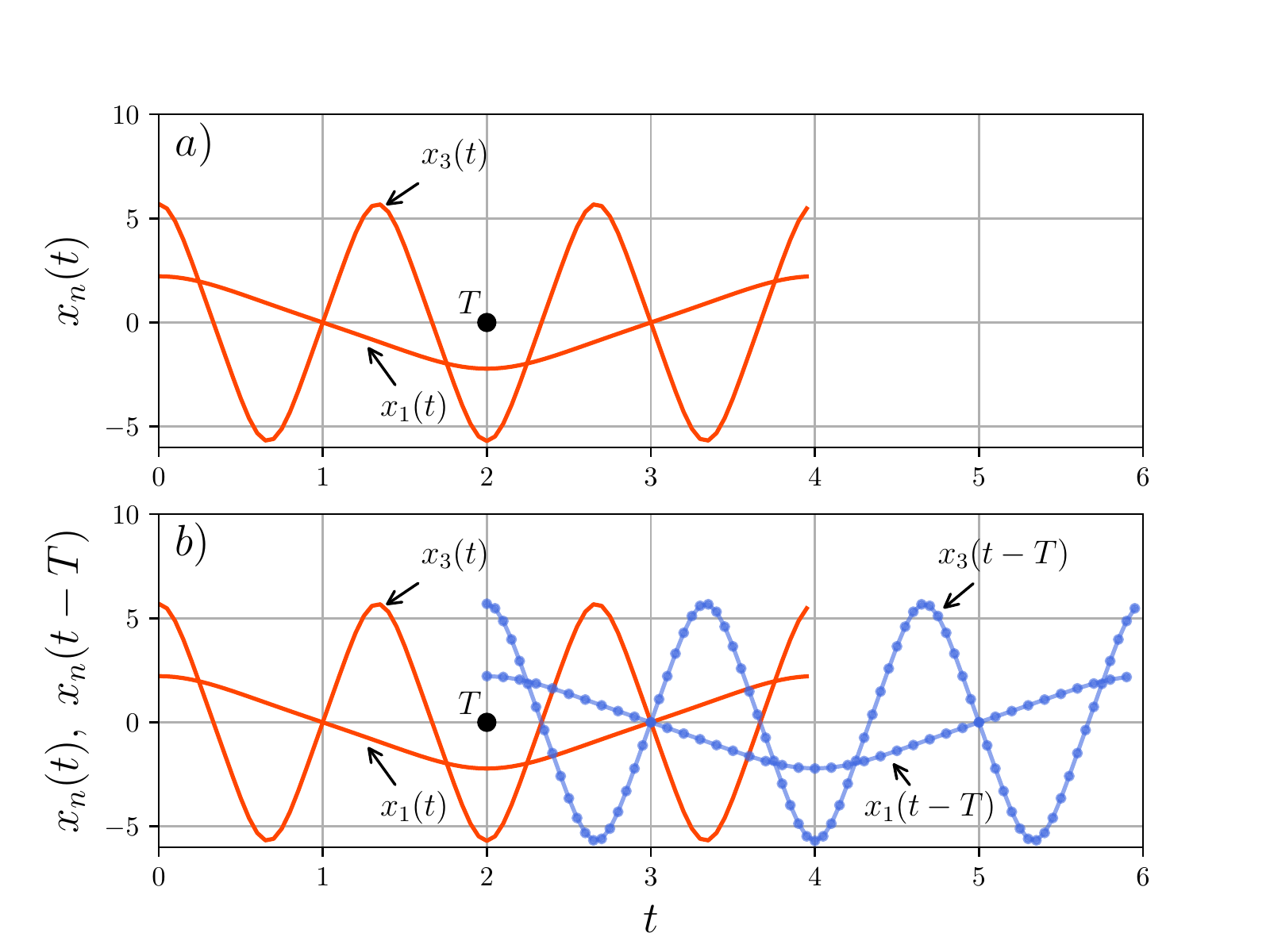} 
	\caption{Solutions of the double-well Duffing ODE (\ref{ode-}), alias odd $n$ in the delayed Duffing DDE (\ref{ode}). Solid red: solutions $x_{n}(t)$ of (\ref{ode-}). Dotted blue: shifted delayed solutions $x_{n}(t-T)$, $T=2$. }
	\label{fig:05}
\end{figure}
%%%%%%%%%%%%%%%%%%%%%%%%%%%%%%

Since each level of positive energy $H>0$ consists of a single periodic orbit $(x,\dot{x})$, with odd force law, the time taken  to travel from any point $(x,\dot{x})$ on a level set to its antipode $(-x,-\dot{x})$ is half its minimal period, $p/2$.
Indeed this fact holds for any odd force law, by time reversibility of the oscillator.
Therefore, every solution of the double-well Duffing ODE (\ref{ode-}) with positive energy $H$ and minimal period $p$ satisfies the oddness symmetry
\begin{equation}
x(t)=-x(t-p/2)\,,
 \label{oddx}
\end{equation}
for all $t$. 

To perform the lift from the double-well Duffing ODE \eqref{ode-} to the delayed Duffing DDE \eqref{dde}, we now fix any time delay $T>0$ (black dot $T=2$ in Figure \ref{fig:05}), this time without any further constraint.
For $p_1:=2T$, the delay $T$ coincides with half the minimal period of the solution $x_1(t)$ of the non-delayed double-well Duffing ODE \eqref{ode-}.
The oddness symmetry \eqref{oddx} at half period $p_1/2=T$ therefore implies that $x_1(t)$ also solves our original delayed Duffing DDE \eqref{dde}, 
\begin{equation}
\ddot{x}+x(t-p/2)+{{x}^{3}}=0 \label{x4} \,.
\end{equation}

Analogously, we can perform the lift from the non-delayed double-well Duffing ODE \eqref{ode-} to the delayed Duffing DDE \eqref{dde}, for any odd $n=1,3,5,\dots$, as follows.
Let $x_n$ denotes the ODE solution of \eqref{ode-} with minimal period $p_n:=2T/n$. 
Then oddness symmetry \eqref{oddx} implies 
\begin{equation}
x_{n}(t-T) = -x_{n}(t),~~~~~~~~~~ \textrm{for all odd }n.
\label{j_T}
\end{equation}
Substitution into \eqref{ode-} implies that $x_n(t)$ also solves \eqref{dde}.
See Figure \ref{fig:05}a,b for illustrations of the cases $n=1,3$.
Note how $p_n\searrow 0$ implies unbounded amplitudes $\sqrt{2}<A_n \nearrow \infty$, for $n\rightarrow \infty$.
Thus we obtain an unbounded sequence of more and more rapidly oscillating periodic solutions to the delayed Duffing DDE \eqref{dde}, with minimal periods $2T, 2T/3, 2T/5, \dots$, which are also periodic with (non-minimal) period $2T$.
This proves our claim \eqref{p_n}, for odd $n$. 

By (\ref{sol}), the exact periodic solutions $x_n(t)$ are expressed as Jacobi elliptic functions
\begin{equation*}
x_{n}(t) = A_{n}\, \mathrm{cn}(\omega_{n}\,t, m_{n}),
\end{equation*}
where (\ref{init_4}) becomes
\begin{equation}
m_{n} = \frac{A_{n}^2}{2(A_{n}^2-1)}~~~~~~~~~\textrm{and} ~~~~~~~~~~ \omega_{n} = \sqrt{A_{n}^2-1}\,.
\label{init_2}
\end{equation}
As we have mentioned in subsection \ref{DuffingODE}, the positivity and symmetry condition $H>0$ becomes equivalent to $A>\sqrt{2}$.  

\begin{figure}[t!]
	\begin{center}
		\def\svgwidth{.8\textwidth}
		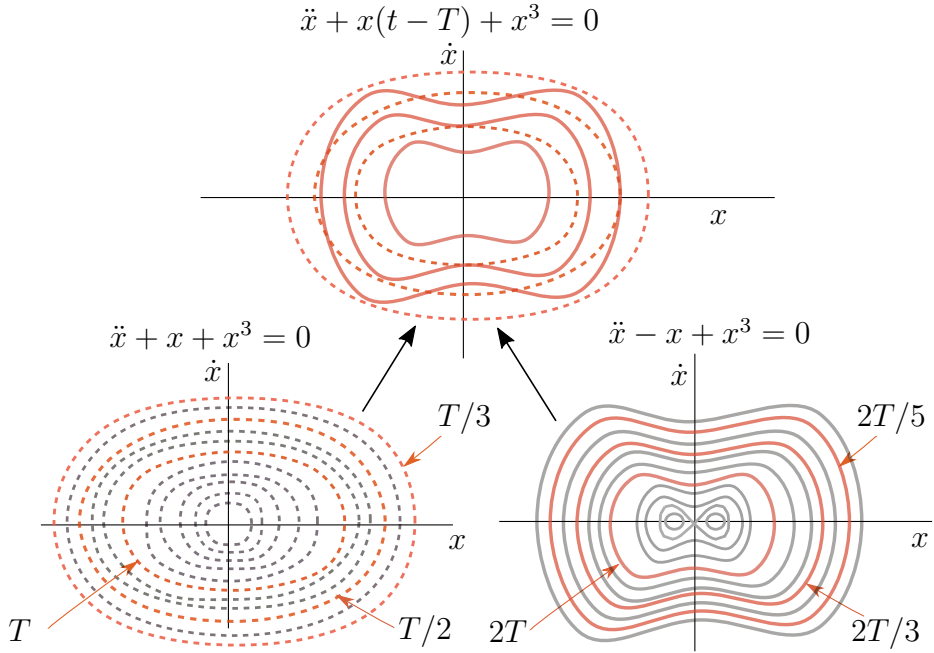   \hspace{0.1cm}
		\caption{Schematic illustration of the lifts from the non-delayed Duffing ODEs \eqref{ode+}, \eqref{ode-} (bottom) to the delayed Duffing DDE (\ref{dde}) (top). }
		  \label{fig:06}
	\end{center}
\end{figure}

Figure \ref{fig:06} schematically illustrate the lift from the non-delayed Duffing ODE (\ref{ode}) to the delayed Duffing DDE (\ref{dde}), for both even and odd $n$. This lift will be used in the next section to numerically determine the amplitudes $A_n$ of the lifted, rapidly oscillating periodic solutions of the DDE (\ref{dde}).

%%%%%%%%%%%%%%%%%%%%%%%%%%%%%
%%%%%%%%%%%%%%%%%%%%%%%%%%%%%

%%%%%%%%%%%%%%%%%%%%%%%%%%%%%%%%%%%%%%%%%%%%%%%%%%%%%%
\section{Amplitudes}\label{ampl}

We sketch two practical approaches to determine the amplitudes $A_n$ of the rapidly oscillating periodic solutions $x_n(t)$ in the delayed Duffing equation (\ref{dde}). 
One approach is essentially numerical; the other approach is analytic, based on an exact series expansion at $n=\infty$ and at infinite amplitude. 
%We discuss the direct numerical approach, first.

The amplitudes $A_n$ of the lifted solutions $x_n(t)$ arise from the closed curves $H>0$ in the non-delayed Duffing ODE (\ref{ode}) with specific values 
\begin{equation}
p \equiv p_n = 2T/n,~~~~~~~~~~~~~~~ n = 1, 2, 3, \dots,
\label{p_2}
\end{equation}
of their minimal period.
See \eqref{p_n} and section \ref{lift} for details.

Substitution of \eqref{p_2} into the explicit elliptic integral \eqref{p_1} provides the implicit equation
\begin{equation}
\label{A_nK}
2T/n = p = 4\,K(m(A_n))/ \omega(A_n)\,
\end{equation}
for $A_n$, given $T$ and $n$.
Here the functions $m(A_n)$ and $\omega(A_n)$ are specified in \eqref{init_4}; we have suppressed explicit dependence on the parity of $n$ in this abbreviated notation.

For high precision numerical solutions $A_n$ of \eqref{A_nK} we rely on the Python-based Newton solver \texttt{fsolve}.
The Newton-method requires initial approximations for the desired solution $A_n$; for initial guesses we use the formal expansions in \cite{Davidow}, Eq.\,(4). The complete elliptic integral $K(m)$ in \eqref{A_nK} is evaluated using the Python-based quadrature \texttt{quad}. The integration is performed using a Clenshaw-Curtis method which uses Chebyshev moments. For $T=3$, for example, the reference amplitudes $A_{n}$ corresponding to the red marked points in Figure \ref{fig:03} are found to be $A_1=1.74566491$\dots,  $A_{2}=2.16089536$\dots, $A_{3}=3.90053028$\dots, and $A_{4}=4.79499435$\dots .

Note that time delays $T$ and $\bar{T}$ share the same reference amplitudes, if the relation  $T/n=\bar{T}/\bar{n}$ holds.
Here $n$ and $\bar{n}$ are required to be both odd, or both even.
For example the amplitude $A_n=A_n(T)$, for $n=1$ and $T=0.1$, coincides with the amplitude $A_{\bar{n}}(\bar{T})$, for $\bar{n}=3$ and $\bar{T}=0.3$. 

Our second approach is analytic in nature.
We start from \eqref{A_nK} with an exact Taylor expansion of $p(A):=4K(m(A))/ \omega(A)$, at $A=\infty$, with respect to $1/A$.
For even $n$ the functions $m(A)$ and $\omega(A)$ have been specified in \eqref{init_1}.
Up to errors of order 13 in $1/A$ we obtain
\begin{align}
\label{p(A)even}
\begin{split}
p = \frac{\gamma}{\sqrt{\pi}} \Big( &
 A^{-1} - 
   \left(\tfrac{1}{2} + 4 \pi^2/\gamma^2 \right) A^{-3} +
   \left(\tfrac{1}{2} + 6 \pi^2/\gamma^2 \right) A^{-5}  -
   \left(\tfrac{5}{8} + 9 \pi^2/\gamma^2 \right) A^{-7} + \\
   &+ \left(\tfrac{85}{96}+ 14 \pi^2/\gamma^2 \right) A^{-9} -
         \left(\tfrac{87}{64} + \tfrac{903}{40} \pi^2/\gamma^2 \right) A^{-11}\Big) +
   \mathcal{O}\left(A^{-13}\right) \,.
\end{split}
\end{align}
Here $\gamma := \Gamma(1/4)^2$ denotes the square of the Euler Gamma-function, evaluated at 1/4.
Note $p=0$ at $A=\infty$. 
Inverting the above series provides an expansion of the inverse function $A(p)$.
Specifically, the Taylor expansion  of $A$ as a function of $1/p$ at $p=0$, up to errors of order 11 in $1/p$ reads
\begin{align}
\label{A(p)even}
\begin{split}
A = \frac{\gamma}{\sqrt{\pi}} \Big(&
 p^{-1} - 
   \pi \left(\tfrac{1}{2}\gamma^2 + 4 \pi^2 \right)\gamma^{-4} p -
   2\pi^4 \left(\gamma^2 + 16 \pi^2 \right)\gamma^{-8} p^3  -\\
   &-8\pi^7\left(3\gamma^2 + 56 \pi^2 \right) \gamma^{-12} p^5 + 
    \tfrac{1}{96} \pi^4 \left(\gamma^8 - 36\,864\,\gamma^2\pi^6 - 737\,280\, \pi^8 \right) \gamma^{-16} p^7 +\\
    &+\tfrac{1}{960} \pi^5 \left(5\gamma^{10} + 328 \gamma^8\pi^2-6\,758\,400\,\gamma^2\pi^8-140\,574\,720\,\pi^{10} \right) \gamma^{-20} p^9 \Big) + \\
   &+\mathcal{O}\left(p^{11}\right) \,.
\end{split}
\end{align}
Inserting $p=2T/n$ readily provides Taylor expansions of $A$ with respect to $n$, in the limit of large $n\rightarrow \infty$ and for any fixed delay $T>0$.
Alternatively, of course, we may consider $n$ fixed and read \eqref{A(p)even} as an expansion with respect to small delays $T>0$, or with respect to any small combination of $T/n$.

For odd $n$, the analogous expansions have to be based on the functions $m(A)$ and $\omega(A)$ specified in \eqref{init_2}.
With the same notation as above we obtain

\begin{align}
\label{p(A)odd}
\begin{split}
p = \frac{\gamma}{\sqrt{\pi}} \Big(&
 A^{-1} + 
   \left(\tfrac{1}{2} + 4 \pi^2/\gamma^2 \right) A^{-3} +
   \left(\tfrac{1}{2} + 6 \pi^2/\gamma^2 \right) A^{-5}  +
   \left(\tfrac{5}{8} + 9 \pi^2/\gamma^2 \right) A^{-7} + \\
   &+ \left(\tfrac{85}{96}+ 14 \pi^2/\gamma^2 \right) A^{-9} +
         \left(\tfrac{87}{64} + \tfrac{903}{40} \pi^2/\gamma^2 \right) A^{-11}\Big) +
   \mathcal{O}\left(A^{-13}\right) \,.
\end{split}
\end{align}

\begin{align}
\label{A(p)odd}
\begin{split}
A = \frac{\gamma}{\sqrt{\pi}} \Big(&
 p^{-1} + 
   \pi \left(\tfrac{1}{2}\gamma^2 + 4 \pi^2 \right)\gamma^{-4} p -
   2\pi^4 \left(\gamma^2 + 16 \pi^2 \right)\gamma^{-8} p^3  +\\
   &+8\pi^7\left(3\gamma^2 + 56 \pi^2 \right) \gamma^{-12} p^5 + 
    \tfrac{1}{96} \pi^4 \left(\gamma^8 - 36\,864\,\gamma^2\pi^6 - 737\,280\, \pi^8 \right) \gamma^{-16} p^7 -\\
    &-\tfrac{1}{960} \pi^5 \left(5\gamma^{10} + 328 \gamma^8\pi^2-6\,758\,400\,\gamma^2\pi^8-140\,574\,720\,\pi^{10} \right) \gamma^{-20} p^9 \Big) + \\
   &+\mathcal{O}\left(p^{11}\right) \,.
\end{split}
\end{align}

Comparing the even and odd cases, we observe how their sign patterns are related by the complex linear transformation $p\mapsto \mathrm{i}p,\ A \mapsto \mathrm{i}A$.
This is in agreement with a scaling of the Duffing ODE.

We emphasize that all Taylor expansions \eqref{p(A)even}--\eqref{A(p)odd} are convergent and hence can be performed up to any order.
Worries like secular terms and other nuisances ubiquitous in formal asymptotics, disappear.
In summary, analytic expansions work best for small $T/n$, e.g. for large $n$, where numerical methods face increasing difficulties.
The numerical approach, on the other hand, is the method of choice for larger $T/n$, e.g. for small $n$.

%%%%%%%%%%%%%%%%%%%%%%%%%%%%%%%%%%%%%%%%%%%%%%%%%%%%%%
\section{Stability}\label{stab}

In this section we summarize results from \cite{Fieetal} on local asymptotic stability and instability of the rapidly oscillating periodic solutions $x_n(t), \ n=1,2,3,\dots$, of the delayed Duffing DDE \eqref{dde}, as constructed in section \ref{lift}.
We recall how the ODE solutions $x_n$ of \eqref{ode} with positive energy $H$ are uniquely determined by their minimal periods $p_n=2T/n$, where $T>0$ denotes the delay in \eqref{dde}; see \eqref{p_n} and \eqref{sol}--\eqref{p_1}.

To be precise we recall that a periodic reference orbit $x_*$ is called \emph{stable limit cycle}, or also \emph{locally asymptotically stable}, if any other solution $x(t)$, which starts sufficiently nearby, remains near the set $x_*$ and converges to that set, for $t \rightarrow\infty$.
A sufficient (but not necessary) condition for local asymptotic stability is \emph{linear asymptotic stability}.
In other words, all \emph{Floquet} (alias \emph{Lyapunov}) \emph{exponents} $\eta$ of the periodic orbit $x_*$ possess strictly negative real part (except for the algebraically simple trivial exponent $\eta=0$).
We speak of \emph{linear instability}, in contrast, if $x_*$ possesses any Floquet (alias Lyapunov) exponent with strictly positive real part.
Deeper results on unstable manifolds then imply \emph{nonlinear instability}.
In fact, there exists a solution $x(t)$ which is defined for all $t\leq 0$ and converges to $x_*$ in backwards time $t\rightarrow -\infty$.

The stability results of \cite{Fieetal} specialize to our present context as follows.

\begin{thm} \label{thmodd}
Let $n$ be odd and assume
\begin{equation}
0<T^2<\tfrac{3}{2}\pi^2.	
\label{Ttorus}
\end{equation}
Moreover assume that $n\geq n_0(T)$ is chosen large enough.

Then the periodic orbit $x_n$ of the delayed Duffing equation \eqref{ode} is asymptotically stable, both linearly and locally.
\end{thm}

\begin{thm} \label{thmeven}
Let $n$ be even, $T>0$, and assume $n\geq n_0(T)$ is chosen large enough.

Then the periodic orbit $x_n$ of the delayed Duffing equation \eqref{ode} is linearly and nonlinearly unstable.
\end{thm}

For the leading Floquet exponent $\eta$, i.e. the nontrivial exponent with real part closest to zero, the precise asymptotics
\begin{equation}
\label{Floq}
\eta = \tfrac{2}{3} (-1)^{n+1}T^2 + \dots
\end{equation}
has been derived, for even and odd $n \rightarrow \infty$.

Towards the stability boundary $T^2=\tfrac{3}{2}\pi^2$ of Theorem \ref{thmodd}, the periodic orbits $x_n$ with odd $n$ lose stability, and undergo a torus bifurcation of Neimark-Sacker-Sell type.
In particular, rational rotation numbers on the bifurcating torus will indicate periodic orbits of the delayed Duffing DDE \eqref{ode} which are \emph{not} lifts of the ODE Duffing orbits $x_n$ studied in the present paper.

We caution the reader that Floquet theory for delay differential equations is not an entirely trivial matter.
Therefore we only illustrate our stability results in the next section, numerically.
For detailed mathematical proofs we have to refer to \cite{Fieetal}.

%%%%%%%%%%%%%%%%%%%%%%%%%%%%%%%%%%%%%%%%%%%%%%%%%%%%%%
\section{Discussion}\label{disc}

%%%%%%%%%%%%%%%%%%%%%%%%%%%
\begin{figure}[t!]
	\centering
	\vspace*{-1cm}
	\includegraphics[width=\textwidth]{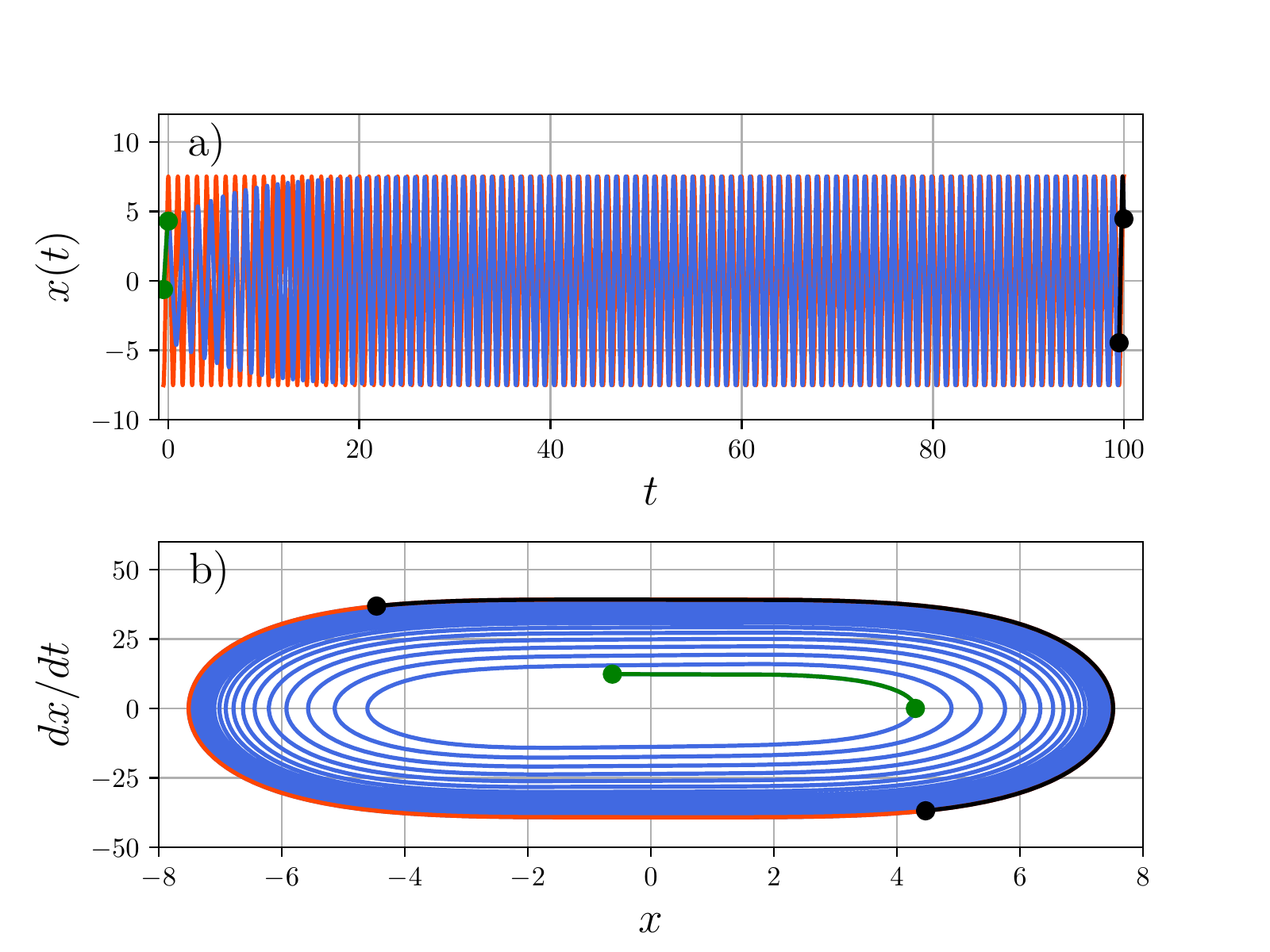} 
	\caption{ Time histories (a) and phase plane plots (b) for delay $T=0.5$. 
	Red: exact periodic solution $x_1(t)$ for $n=1$, with reference amplitude $A_{1}=7.5139958\dots$ and minimal period $2T$; see \eqref{sol}. 
	Blue: simulated solution of the delayed Duffing DDE (\ref{dde}) with initial history function (\ref{init_3}) and initial amplitude $A = 4.3$. 
	Green: initial history function (\ref{init_3}). 
	Black: final state of the history function (\ref{init_3}).  
	Note the convergence of the blue solution to the locally asymptotically stable red limit cycle $x_{1}$, for large times $t$.}
	\label{fig:07}
\end{figure}
%%%%%%%%%%%%%%%%%%%%%%%%%%%%%%

Figure  \ref{fig:07} plots two solutions of  the delayed Duffing equation (\ref{dde}) with delay $T= 0.5$:
a numerical solution $x(t)$ (blue), and the lifted exact solution $x_1(t)$ (red) specified in (\ref{sol}). 
The minimal period $p_1$ of $x_1(t)$ coincides with $2T$; see \eqref{p_n}. 
Figure \ref{fig:07} contains the time history (a) and the phase plane (b).  
The green curve denotes the initial history function
\begin{equation}
\left(x(t), \dot{x}(t)\right) = \left( A\, \mathrm{cn}(\omega\,t,m), -A\,\omega\, \mathrm{sn}(\omega\,t,m)\, \mathrm{dn}(\omega\,t,m)  \right),
\label{init_3}
\end{equation}
for $-T<t<0$ and with initial amplitude $A = 4.3$. 
The values of $m$ and $\omega$ are obtained from (\ref{init_4}) with $n=1$. 
Note how $x(t)$ is a solution of the non-delayed Duffing ODE \eqref{ode} with minimal period $p = p(A) = 1.7972608\dots$.
However, the initial history function $x(t)$ is \emph{not} a solution of the delayed Duffing DDE (\ref{dde}), because $T = 0.5$ is \emph{not} an integer multiple of the larger ODE period $p = 1.7972608\dots$  . 
Therefore the simulated solution $x(t)$ of the  delayed Duffing DDE \eqref{ode} (blue), is \emph{not} periodic.

Instead, the simulated solution (blue), with initial amplitude $A = 4.3$, approaches the exact periodic solution $x_1(t)$ (red) of minimal period $2 T$ and with amplitude $A_1=7.5139958\dots$ . 
Indeed, the black curve indicates the history function, for $100-T<t<100$, of the final state of the blue solution $x(t)$ at $t=100$.
The stability result of Theorem \ref{thmodd} only asserts local convergence to $x_n$ for large odd $n$, but not for $n=1$.
The convergence to $x_1$ indicates how that stability result might actually extend, all the way, down to  the smallest possible choice $n=1$.
Moreover, ``local'' attraction to $x_1$ holds sway over quite a distance, down to an initial amplitude $A=4.3$ significantly smaller than the asymptotic amplitude $A_1=7.5139958\dots$ of $x_1$.

%%%%%%%%%%%%%%%%%%%%%%%%%%%
\begin{figure}[t!]
	\centering
	\vspace*{-1cm}
	\includegraphics[width=\textwidth]{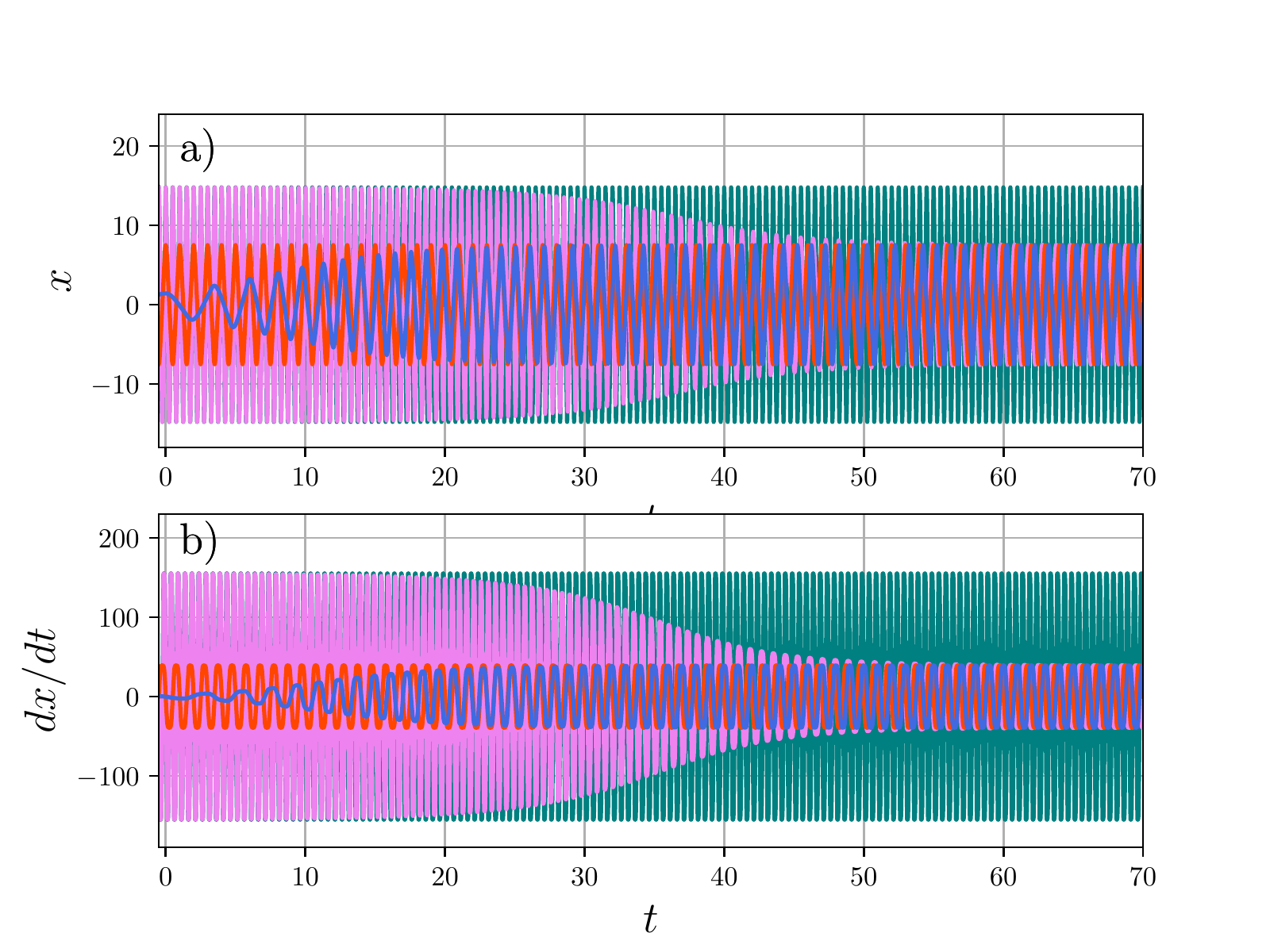}
	\caption{ Time histories of $x(t)$, top (a), and of $\dot{x}(t)$, bottom (b), for delay $T=0.5$. 
	Exact solutions $x_n(t)$ for $n=1$ (red) and $n=2$ (teal); see (\ref{sol}). 
	Their amplitudes are $A_1=7.5139958\dots$ and $A_2 = 14.7834172\dots$, respectively. 
	The numerical solution of the delayed Duffing DDE (\ref{dde}) with initial amplitude $A= 1.42$ (blue) illustrates wide asymptotic stability of the stable limit cycle $x_1$.
	The numerical solution with initial amplitude $A= 14.77$ (violet), quite close to $A_2$, indicates a heteroclinic orbit from the unstable periodic orbit $x_2$ to the stable limit cycle $x_1$.}
	\label{fig:08}
\end{figure}
%%%%%%%%%%%%%%%%%%%%%%%%%%%%%%

Figure \ref{fig:08}  compares two lifted exact periodic solutions, $x_1(t)$ (red) and $x_2(t)$ (teal).
Two numerical solutions of the delayed Duffing DDE (\ref{dde}) for $T=0.5$ are included, which arise from the two initial history functions (\ref{init_3}) with initial amplitudes $A= 1.42$ (blue) and $A= 14.77$ (violet), respectively. 
The reference amplitudes corresponding to the exact $n=1$ (red) and $n=2$ (teal) periodic solutions (\ref{sol}) are $A_1=7.5139958\dots$ and $A_2 = 14.7834172\dots$, respectively. 
Figure \ref{fig:08} indicates how both simulated solutions (blue and violet) approach the exact stable limit cycle $x_1$ (red); see Theorem \ref{thmodd}. 
Also note how the simulated solution with initial condition $A=14.77$ (violet) starts very close to the exact, but linearly unstable, periodic solution (teal) of $A_2 = 14.7834172\dots$, but eventually diverges as time $t$ increases.
See Theorem \ref{thmeven}.
This indicates the presence of a heteroclinic orbit $x(t)$, from $x_2$ to $x_1$, which is defined for all positive and negative times $t$ and converges to $x_2$, for decreasing $t \searrow -\infty$, and to $x_1$, for increasing $t \nearrow +\infty$. 

Our periodicity Ansatz requires half minimal periods $p/2=T/n$ to be integer fractions $n=1,2,3, \dots$ of the delay $T$. 
Of course we have to caution the reader that there may be many periodic solutions of the DDE (\ref{dde})  which are not captured by this Ansatz.

%%%%%%%%%%%%%%%%%%%%%%%%%%%%%%%%%%%%%%%%%%%%%%%%%%%%%%
\section{Conclusion}\label{conc}

In this work we showed how the Duffing equation (\ref{dde}) with time delay $T$ possesses an unbounded sequence of infinitely many rapidly oscillating periodic solutions $x_n(t),\ n=1,2,3,\dots$ . 

Each solution $x_n$ arises from a periodic solution $x_n(t)$ of the non-delayed classical Duffing equation (\ref{ode}) with minimal period $p_n=2T/n$.  
In particular, the classical non-delayed Duffing oscillator provides an unbounded sequence of exact periodic solutions of the delayed Duffing equation.  
Based on the Hamiltonian energy of the classical Duffing equation, and standard Jacobi elliptic integrals, we have also derived high-precision reference amplitudes of these periodic solutions $x_n$.

For delays $T$ such that $0<T^2<\tfrac{3}{2} \pi^2$, and for odd $n$ large enough, the solutions $x_n$ are locally asymptotically stable limit cycles.
For large even $n$, in contrast, the solutions $x_n$ are linearly and nonlinearly unstable.

We have illustrated our results with numerical simulations, for low $n=1,2$.

%%%%%%%%%%%%%%%%%%%%%%%%%%%%%%%%%%%%%%%%%%%%%%%%%%%%%%

\end{document}